\documentclass[12pt]{article}
\title{A strong failure of $\aleph_0$-stability for atomic classes}
\author{Michael C.\
Laskowski\thanks{Partially supported
by NSF grant DMS-1308546.}\\
Department of Mathematics\\University of Maryland
\and
Saharon Shelah\thanks{
Partially supported by European Research Council grant 338821 and
NSF grant DMS-1362974.
 Publication no.\ 1099.}
 \\Hebrew University\\ Rutgers University
}

\usepackage{amssymb}
\usepackage{url}

\def\abar{\bar{a}}
\def\bbar{\bar{b}}
\def\cbar{\bar{c}}
\def\dbar{\bar{d}}
\def\ebar{\bar{e}}
\def\mbar{\overline{m}}
\def\Mbar{\overline{M}}
\def\Nbar{\overline{N}}
\def\pbar{\bar{p}}
\def\qbar{\bar{q}}
\def\rbar{\bar{r}}

\def\xbar{\overline{x}}

\def\phi{\varphi}

\def\B{{\cal B}}

\def\F{{\cal F}}
\def\FF{{\bf F}}

\def\P{{\cal P}}
\def\PP{{\mathbb P}}
\def\Q{{\mathbb Q}}

\def\S{{\cal S}}

\def\tp{{\rm tp}}

\def\dom{{\rm dom}}

\def\pcl{{\rm pcl}}

\def\Fa0{{\FF^a_{\aleph_0}}}

\def\bp{{\bf Proof.}\quad}
\def\endproof{\medskip}
\def\<{\langle}
\def\>{\rangle}
\def\o2{{^{\omega} 2}}

\def\n2{{^{n} 2}}

\newtheorem{Theorem}{Theorem}[section]
\newtheorem{Proposition}[Theorem]{Proposition}
\newtheorem{Definition}[Theorem]{Definition}
\newtheorem{Notation}[Theorem]{Notation}

\newtheorem{Lemma}[Theorem]{Lemma}
\newtheorem{Corollary}[Theorem]{Corollary}

\newtheorem{Conclusion}[Theorem]{Conclusion}

\def\At{{\bf At}}
\def\rk{{\rm rk}}
\def\QQ{{\cal Q}}
\def\CC{{\cal C}}
\def\SS{{\cal S}}
\def\pp{{\bf p}}
\def\qq{{\bf q}}
\def\rr{{\bf r}}
\def\ss{{\bf s}}
\def\aabar{{\bf \overline{a}}}
\def\nbar{{\overline{n}}}
\def\U{{\cal U}}
\def\spl{{\rm spl}}
\def\tr{{\rm tr}}

\newcommand\myrestriction{\mathord\restriction}
\def\mr#1{\myrestriction_{#1}}

\date{\today}

\begin{document}

\maketitle

\begin{abstract}  We study classes of atomic models  $\At_T$ of a countable, complete first-order theory $T$.  We prove that if
$\At_T$  is not $\pcl$-small, i.e.,  there is an atomic  model $N$ that realizes uncountably many types over $\pcl(\abar)$ for some finite $\abar$ from $N$,
then there are $2^{\aleph_1}$ non-isomorphic atomic models of $T$, each of size $\aleph_1$.
\end{abstract}

\section{Introduction}\label{intro}

In  a series of papers \cite{BLS,BL,BLS3}, Baldwin and the authors have begun to develop a model theory for complete sentences of $L_{\omega_1,\omega}$ that have fewer than $2^{\aleph_1}$ non-isomorphic models of size $\aleph_1$.  By well known reductions, one can replace the  reference to infinitary sentences by restricting to the class of {\em atomic} models of a countable, complete first-order theory.\footnote{Specifically, for every complete sentence $\Phi$ of $L_{\omega_1,\omega}$, there is a complete first-order theory $T$ in a countable vocabulary containing the vocabulary of $\Phi$ such that the models of $\Phi$ are precisely the reducts of the class of atomic models of $T$ to the smaller vocabulary.}

Fix, for the whole of this paper, a complete theory $T$ in a countable language that has at least one atomic model\footnote{A model $M$ is {\em atomic} if, for every finite
tuple $\abar$ from $M$, $\tp_M(\abar)$ is {\em principal} i.e., is uniquely determined by a single formula $\phi(\xbar)\in\tp_M(\abar)$.} of size $\aleph_1$.  
By theorems of Vaught, 
these restrictions on $T$ are well understood.
Such a $T$ has an atomic model if and only if every consistent formula can be extended to a complete formula.  Furthermore, any two countable, atomic models of $T$ are isomorphic, and a model is prime if and only if it is countable and atomic.  Using a well-known union of chains argument, $T$ has an atomic model of size $\aleph_1$ if and only if the countable atomic model is not minimal, i.e., it has a proper elementary substructure.  

The analysis of $\At_T$, the class of atomic models of $T$, begins by restricting the notion of types to those that can be realized in an atomic model.  
Suppose $M$ is atomic and $A\subseteq M$.  We let
$S_{at}(A)$ denote the set of complete types $p$ over $A$ for which $Ab$ is an atomic set for some (equivalently, for every) realization $b$ of $p$.  It is easily checked that when $A$ is countable, $S_{at}(A)$ is a $G_\delta$ subset of the Stone space $S(A)$, hence $S_{at}(A)$ is Polish with respect to the induced topology.    We will repeatedly use the fact that any countable, atomic set $A$ is contained in a countable, atomic model $M$.  However, unlike the first-order case, some types in $S_{at}(A)$ need not extend to types in $S_{at}(M)$.
Indeed, there are examples where the space $S_{at}(A)$ is uncountable (hence contains a perfect set) while $S_{at}(M)$ is countable.  
Thus, for analyzing types over countable, atomic sets $A\subseteq M$, we are led to consider 
$$S_{at}^+(A,M):=\{p|A:p\in S_{at}(M)\}.$$

Equivalently, $S_{at}^+(A,M)$ is the set of $q\in S_{at}(A)$ that can be extended to a type $q^*\in S_{at}(M)$.

Next, we recall the notion of {pseudo-algebraicity}, which was introduced in \cite{BLS}, that is the correct analog of algebraicity in the context of atomic models.  Suppose $M$ is an atomic model, and $b,\abar$ are from $M$.  We say $b\in\pcl_M(\abar)$  if $b\in N$ for every elementary submodel $N\preceq M$ that contains $\abar$.  The seeming dependence on $M$ is illusory -- as is noted in \cite{BLS}, if $b',\abar'$ are inside another atomic model $M'$,  and $\tp_{M'}(b'\abar')=\tp_M(b\abar)$, then $b\in\pcl_M(\abar)$  if and only if
$b'\in\pcl_{M'}(\abar')$.  It is easily seen that inside any atomic model $M$, $\pcl_M(\abar)$ is countable for any finite tuple $\abar$.  Moreover, if $f:M\rightarrow M'$ is an isomorphism of atomic models, then $f(\pcl_M(\abar))=\pcl_{M'}(f(\abar))$ setwise.  As an important special case, if  $\abar\subseteq M'\preceq M$ and  $f:M\rightarrow M'$ fixes $\abar$ pointwise, then
$f$ induces an elementary permutation on $D=\pcl_{M}(\abar)$, which in turn induces a bijection between $S_{at}^+(D,M)$ and $S_{at}^+(D,M')$.

We now give the major new definition of this paper:

\begin{Definition}  {\em An atomic class $\At_T$ with an uncountable model is {\em pcl-small} if, for every atomic model $N$ and for every finite $\abar$ from $N$,
$N$ realizes only countably many complete types over $\pcl_N(\abar)$.
}
\end{Definition}

The name of this notion is by analogy with the first-order case -- A complete, first-order theory $T$ is small if and only if for every model $N$ and every finite $\abar$ from $N$,
$N$ realizes only countably many complete types over $\abar$.
%
%
%
The following proposition relates $\pcl$-smallness with the spaces of types $S^+_{at}(D,M)$.

\begin{Proposition}  \label{relate}  The atomic class $\At_T$ is $\pcl$-small if and only if the space of types $S^+_{at}(\pcl_M(\abar),M)$ is countable for every countable, atomic model $M$ and every finite $\abar$ from $M$.
\end{Proposition}

\bp  
First, assume that some atomic model $N$ and finite sequence $\abar$ from $N$ witness that $\At_T$ is not $\pcl$-small.  Choose $\{c_i:i\in\omega_1\}\subseteq N$
realizing distinct complete types over $D=\pcl_N(\abar)$.  Also, choose a countable $M\preceq N$ that contains $\abar$, and hence $D$.  Then $\{\tp(c_i/D):i\in\omega_1\}$ witness that 
$S^+_{at}(D,M)$ is uncountable.

For the converse, choose a countable, atomic model $M$ and $\abar$ from $M$ such that $S^+_{at}(D,M)$ is uncountable, where $D=\pcl_M(\abar)$.  
We will inductively construct a continuous, increasing elementary chain $\<M_\alpha:\alpha<\omega_1\>$ of countable, atomic models with $M=M_0$ and, for each ordinal $\alpha$, there is an element $c_\alpha\in M_{\alpha+1}$
such that $\tp(c_\alpha/D)$ is not realized in $M_\alpha$.  Given such a sequence, it is evident that $N=\bigcup_{\alpha<\omega_1} M_\alpha$ and $\abar$ witness that $\At_T$ is not $\pcl$-small.  To construct such a sequence, we have defined $M_0$ to be $M$ and take unions at limit ordinals.  For the successor step, assume $M_\alpha$ has been defined.
As $M$ and $M_\alpha$ are each countable atomic models that contain $\abar$, choose an isomorphism $f:M\rightarrow M_\alpha$ fixing $\abar$ pointwise.  As noted above,
$f$ fixes $D$ setwise.  As $M_\alpha$ is countable, so is the set $\{\tp(c/D):c\in M_\alpha\}$.  As $S^+_{at}(D,M)$ is uncountable, choose an atomic type
$p\in S_{at}(M)$, whose restriction to $D$  is distinct from $\{f^{-1}(\tp(c/D)):c\in M_\alpha\}$.  Now choose $c_\alpha$ to realize $f(p)$.  Then, as $M_\alpha c_\alpha$
is a countable atomic set, choose a countable elementary extension $M_{\alpha+1}\succeq M_\alpha$ containing $c_\alpha$.  
\endproof

Recall that an atomic class $\At_T$ is $\aleph_0$-stable\footnote{Sadly, this usage of `$\aleph_0$-stability' is analogous, but distinct from, the familiar first-order notion.}
 if $S_{at}(M)$ is countable for all (equivalently, for some) countable atomic models $M$.  
As $S^+_{at}(A,M)$ is a set of projections of types in $S_{at}(M)$,  it will be countable whenever $S_{at}(M)$ is.  This observation makes the  following corollary  to  
Proposition~\ref{relate} immediate:

\begin{Corollary} \label{omega} If an atomic class $\At_T$ is $\aleph_0$-stable, then $\At_T$ is $\pcl$-small.
\end{Corollary}

The converse to Corollary~\ref{omega} fails.  For example, the theory $T=REF(bin)$ of countably many, binary splitting equivalence relations is not $\aleph_0$-stable, yet
$\pcl_M(\abar)=\abar$ for every model $M$ and $\abar$ from $M$.  Thus, $S_{at}(\pcl_{M}(\abar))$ and hence  $S^+_{at}(\pcl(\abar),M)$ is countable for every finite tuple $\abar$ inside any atomic model $M$.
The main theorem of this paper is:

\begin{Theorem}  \label{big}  Let $T$ be a countable, complete theory $T$ with an uncountable atomic model.  If the atomic class $\At_T$ is not $\pcl$-small, then
there are $2^{\aleph_1}$ non-isomorphic models in $\At_T$, each of size $\aleph_1$.
\end{Theorem}

Section~2 sets the stage for the proof.  It describes the spaces of types $S^+_{at}(A,M)$, states a transfer theorem for sentences of $L_{\omega_1,\omega}(Q)$, and
details a non-structural configuration arising from non-$\pcl$-smallness.  In Section~3, the non-structural configuration is exploited to give a family of $2^{\aleph_0}$
non-isomorphic structures $(N,\bbar^*)$, where each of the reducts $N$ is in $\At_T$ and has size $\aleph_1$.
Theorem~\ref{big} is finally proved in Section~4.  It is remarkable that whereas it is a ZFC theorem, the proof is non-uniform depending on the relative sizes of
the cardinals
$2^{\aleph_0}$ and $2^{\aleph_1}$. 

\section{Preliminaries}

In this section, we develop some general tools that will be used in the proof of Theorem~\ref{big}.  

\subsection{On  $S^+_{at}(A,M)$}

In this subsection we explore the space of types 
$$S^+_{at}(A,M)=\{p|A:p\in S_{at}(M)\}$$ where $A$ is a subset of a countable, atomic model $M$.  

Fix a countable, atomic model $M$ and an arbitrary subset $A\subseteq M$.  Let $\P$ denote the space of complete types  in one free variable
over finite subsets of $M$.  As $M$ is atomic, 
$\P$ can be identified with the set of complete formulas $\phi(x,m)$ over $M$.  Implication gives a natural partial order on $\P$, namely
$p\le q$ if and only if $\dom(p)\subseteq\dom(q)$ and $q\vdash p$.   
One should think of elements of $\P$ as `finite approximations' of types in $S_{at}^+(A,M)$.  We describe two conditions on $p\in\P$ that identify extreme behaviors in this regard.

\begin{Definition}  \label{lies}  {\em  
We say a type $p^*\in S^+_{at}(A, M)$ {\em lies above $p\in\P$} if there is some $\pbar\in S_{at}(M)$ extending $p\cup p^*$.  As every $p\in\P$ extends to a type in $S_{at}(M)$,
it follows that at least one $p^*\in S^+_{at}(A,M)$ lies above $p$.
\begin{itemize}
\item  An element $p\in\P$ {\em determines a type in $S^+_{at}(A,M)$} if exactly one $p^*\in S_{at}^+(A,M)$ lies above $p$.
\item  An element $p\in\P$ is {\em $A$-large} if $\{p^*\in S^+_{at}(A,M): p^*$ lies above $p\}$ is uncountable.
\end{itemize}
}
\end{Definition}

To understand these extreme behaviors, we define a rank function $\rk_A:\P\rightarrow (\omega_1+1)$ as follows:
\begin{itemize}
\item  $\rk_A(p)\ge 0$ for all $p\in\P$;
\item  For $\alpha\le\omega_1$, $\rk_A(p)\ge\alpha$ if and only if for every $\beta<\alpha$ and for all finite $F$, $\dom(p)\subseteq F\subseteq M$, there is $q\in S_{at}(F)$
with $q\ge p$
that {\em $\beta$-$A$ splits}, where:
\begin{itemize}
\item  A type $q\in S_{at}(F)$ {\em $A$-splits} if, for some $\phi(x,\abar)$ with $\abar$ from $A$, there are $q_1,q_2\ge q$ with $q\cup\phi(x,\abar)\subseteq q_1$ and
$q\cup \neg\phi(x,\abar)\subseteq q_2$; and $q\in S_{at}(F)$ {\em $\beta$-$A$ splits} if, in addition, $\rk_A(q_1),\rk_A(q_2)\ge\beta$.
\end{itemize}
\item  For $\alpha<\omega_1$, we say $\rk_A(p)=\alpha$ if $\rk_A(p)\ge\alpha$, but $\rk _A(p)\not\ge\alpha+1$.
\end{itemize}

\begin{Proposition}  \label{det}
If $p\in\P$ and $\rk_A(p)=\alpha<\omega_1$, then some $r\ge p$ determines a type in $S^+_{at}(A,M)$.
\end{Proposition}

\bp  We prove this by induction on $\alpha$.  We begin with $\alpha=0$.  Suppose $\rk_A(p)=0$.  As $\rk_A(p)\not\ge 1$, there is a finite $F$, $\dom(p)\subseteq F\subseteq M$
for which there is no $q\in S_{at}(F)$ and $\phi(x,\abar)$ with $\abar$ from $A$ for which $q\ge p$ and both $q\cup\{\phi(x,\abar)\}$ and $q\cup\{\neg\phi(x,\abar)\}$ are consistent.  
So fix any $r\in S_{at}(F)$ with $r\ge p$.  Any such $r$ determines a type in $S^+_{at}(A,M)$.

Next, choose $0<\alpha<\omega_1$ and assume the Proposition holds for all $\beta<\alpha$.  Choose $p\in S_{at}(E)$ with $\rk_A(p)=\alpha$.  As $\rk_A(p)\ge\alpha$, while $\rk_A(p)
\not\ge \alpha+1$, there is a finite $F$, $E \subseteq F\subseteq M$ for which there is no $q\in S_{at}(F)$ that both extends $p$ and $\alpha$-$A$ splits.  So choose any $q\in S_{at}(F)$ with $q\ge p$.  If $q$ determines a type in $S^+_{at}(A,M)$, then we finish, so assume otherwise.  Thus, there is some $\phi(x,\abar)$ with $\abar$ from $A$ such that both $q\cup\{\phi(x,\abar)\}$ and $q\cup\{\neg\phi(x,\abar)\}$ are consistent.  Choose complete types $q_1,q_2\in S_{at}(F\abar)$ extending these partial types.  Clearly, both $q_1,q_2\ge q$,
but since $q$ does not $\alpha$-$A$ split, at least one of them has $\rk_A(q_{\ell})<\alpha$.  But then by our inductive hypothesis, there is $r\ge q_{\ell}$ that determines a type
in $S^+_{at}(A,M)$ and we finish.
\endproof

Next, we turn our attention to $A$-large types and types of rank at least $\omega_1$ and see that these coincide.
We begin with two lemmas, the first involving types of rank at least $\omega_1$ and the second involving $A$-large types.

\begin{Lemma}  \label{presrk}
Assume  that $E\subseteq M$ is finite and $p\in S_{at}(E)$ has $\rk_A(p)\ge\omega_1$.  Then:
\begin{enumerate}
\item  For every finite $F$, $E\subseteq F\subseteq M$, there is $q\in S_{at}(F)$, $q\ge p$, with $\rk_A(q)\ge\omega_1$; and
\item  There is some formula $\phi(x,\abar)$ with $\abar$ from $A$ and $q_1,q_2\in\P$ with $p\cup \{\phi(x,\abar)\}\subseteq q_1$, $p\cup \{\neg\phi(x,\abar)\}\subseteq q_2$,
and both $\rk_A(q_1),\rk_A(q_2)\ge\omega_1$.
\end{enumerate}
\end{Lemma}

\bp (1)  Fix a finite $F$ satisfying $E\subseteq F\subseteq M$.  As $\rk_A(p)\ge\omega_1$, for every $\beta<\omega_1$ there is some $q\ge p$ with $q\in S_{at}(F)$
for which certain extensions of $q$ have rank at least $\beta$.  It follows that $\rk_A(q)\ge\beta$ for any such witness.  However, as $S_{at}(F)$ is countable, there is some
$q\in S_{at}(F)$ which serves as a witness for uncountably many $\beta$.  Thus, $\rk_A(q)\ge\omega_1$ for any such $q\ge p$.

(2)  Assume that there were no such formula $\phi(x,\abar)$.  Then, for any formula $\phi(x,\abar)$, since $\P$ is countable, there would be an ordinal $\beta^*<\omega_1$ such that 
{\bf either}  every $q\in\P$ extending $p\cup\{\phi(x,\abar)\}$, $\rk_A(q)<\beta^*$ {\bf or} every $q\in \P$ extending $p\cup\{\neg\phi(x,\abar)\}$ has $\rk_A(q)<\beta^*$.
Continuing, as there are only countably many formulas $\phi(x,\abar)$, there would be an ordinal $\beta^{**}<\omega_1$ that works for all formulas $\phi(x,\abar)$.
Restating this, $p$ does not $\beta^{**}$-$A$ split, so no extension of $p$ could $\beta^{**}$-$A$ split either.  This contradicts $\rk_A(p)\ge\beta^{**}+1$.
\endproof

\begin{Lemma}
\label{large}  Suppose $q\in S_{at}(F)$ is $A$-large.  Then:
\begin{enumerate}
\item  For every finite $F'$, $F\subseteq F'\subseteq M$, there is some $A$-large $r\in S_{at}(F')$ with $r\ge q$; and
\item  For some $\phi(x,\abar)$, there are $A$-large extensions $r_1\supseteq q\cup\{\phi(x,\abar)\}$ and $r_2\supseteq q
\cup\{\neg\phi(x,\abar)\}$.
\end{enumerate}
\end{Lemma}

\bp Fix such a $q$ and let
 $\S=\{p^*\in S^+_{at}(A,M): p^*$ lies above $q\}$.  
 
 (1) is immediate, since $\S$ is uncountable, while $S_{at}(F')$ is countable.
 
 For (2), first note that 
if there is no such $\phi(x,\abar)$, then there is at most one $p^*\in\S$ with the property that:
\begin{quotation}
\noindent
For any formula $\phi(x,\abar)$ with parameters from $A$, 
$\phi(x,\abar)\in p^*$ if and only if there is an $A$-large $r\in S_{at}(F\abar)$ extending $q\cup\{\phi(x,\abar)\}$.
\end{quotation}
It follows that for any $q^*\in \S-\{p^*\}$, $q^*$ lies over some $r\ge q$ that is not $A$-large.  That is, using the fact that there are only countably many $r\ge q$,
$\S-\{p^*\}$ is contained in the union of countably many countable sets.  But this contradicts $q$ being $A$-large.
\endproof

\begin{Proposition}  \label{largeeq}  For $p\in\P$, $\rk_A(p)\ge\omega_1$ if and only if $p$ is $A$-large.
\end{Proposition}

\bp  First, assume that $\rk_A(p)\ge\omega_1$.  Fix an enumeration $\{c_n:n\in\omega\}$ of $M$.   Using Clauses~(1) and (2) of Lemma~\ref{presrk}, we inductively construct a tree $\{p_\nu:\nu\in 2^{<\omega}\}$ of elements of $\P$ satisfying:
\begin{enumerate}
\item  $\rk_A(p_\nu)\ge\omega_1$ for all $\nu\in 2^{<\omega}$;
\item  If $\lg(\nu)=n$, then $\{c_i:i<n\}\subseteq\dom(p_\nu)$;
\item  $p_{\<\>}=p$;
\item  For $\nu\trianglelefteq\mu$, $p_\nu\le p_\mu$;
\item   For each $\nu$ there is a formula $\phi(x,\abar)$ with $\abar$ from $A$ such that $\phi(x,\abar)\in p_{\nu 0}$ and $\neg\phi(x,\abar)\in p_{\nu 1}$.
\end{enumerate}
Given such a tree, for each $\eta\in 2^\omega$, let $\pbar_\eta:=\bigcup\{p_{\eta |n}:n\in\omega\}$ and let $p^*_\eta:=\pbar_\eta | A$.
By Clauses (2) and (4), each $\pbar_\eta\in S_{at}(M)$, so each $p^*_\eta\in S^+_{at}(A,M)$.   By Clause (5), $p^*_\eta\neq p^*_{\eta'}$ for distinct $\eta,\eta'\in 2^\omega$.
Finally, each of these types lies over $p$ by Clause~(3).  Thus, $p$ is $A$-large.

Conversely, we argue by induction on $\alpha<\omega_1$ that:
\begin{quotation}  $(*)_\alpha:$  \phantom{X} If $p\in \P$ is $A$-large, then $\rk_A(p)\ge\alpha$.
\end{quotation}

Establishing $(*)_0$ is trivial, and for limit $\alpha<\omega_1$, it is easy to establish $(*)_\alpha$ given that $(*)_\beta$ holds for all $\beta<\alpha$.  
So assume $(*)_\alpha$ holds and we will establish $(*)_{\alpha+1}$.
Choose any $A$-large $p\in\P$.  Towards showing $\rk_A(p)\ge\alpha+1$, choose any finite $F$, $\dom(p)\subseteq F\subseteq M$.
As $S_{at}(F)$ is countable and uncountably many types in $S^+_{at}(A,M)$ lie above $p$, there is some $A$-large $q\in S_{at}(F)$ with
$q\ge p$.

Next, by Lemma~\ref{large} choose a formula $\phi(x,\abar)$ with $\abar$ from $A$ such that there are $A$-large extensions $r_1\supseteq q\cup\{\phi(x,\abar)\}$
and $r_2\supseteq q\cup\{\neg\phi(x,\abar)\}$.   Applying $(*)_\alpha$ to both $r_1,r_2$ gives $\rk_A(r_1),\rk_A(r_2)\ge\alpha$.  Thus, $q$ $\alpha$-$A$ splits.
Thus, by definition of the rank, $\rk_A(p)\ge\alpha+1$.
\endproof

We obtain the following Corollary, which is analogous to the statement `If $T$ is small, then the isolated types are dense' from the first-order context.

\begin{Corollary}  \label{countable}  If $S^+_{at}(A,M)$ is countable, then every $p\in\P$ has an extension $q\ge p$ that determines a type in $S^+_{at}(A,M)$.
\end{Corollary}

\bp  If $S^+_{at}(A,M)$ is countable, then no $p\in \P$ is $A$-large.  Thus, every $p\in\P$ has $\rk_A(p)<\omega_1$ by Proposition~\ref{largeeq}, so has an extension determining a type in $S^+_{at}(A,M)$ by Proposition~\ref{det}.
\endproof

We close with a complementary result  about extensions of  $A$-large types.

\begin{Definition}  \label{perfectx}  {\em 
A type $r\in S_{at}(M)$ is {\em $A$-perfect} if $r\mr{A}$ is omitted in $M$ and for every finite $\mbar$ from $M$, the restriction $r\mr{\mbar}$ is $A$-large.
}
\end{Definition}

The name {\em perfect} is chosen because, relative to the usual topology on $S_{at}(M)$, there are a perfect set of 
$A$-perfect types extending any $A$-large $p\in\P$.  However, for what follows, all we need to establish 
is that there are uncountably many, which is notationally simpler to prove.

\begin{Proposition}  \label{perfect}  Suppose $p\in\P$ is $A$-large.  Then there are uncountably many $A$-perfect $r\in S_{at}(M)$ extending $p$.
\end{Proposition}

\bp  Fix an $A$-large $p\in\P$.  Choose a set $R\subseteq S_{at}(M)$ of  representatives for $\{p^*\in S_{at}^+(A,M): p^*$ lies above $p\}$, 
i.e., for every such $p^*$, there is exactly one $\pbar\in R$ whose restriction $\pbar\mr{A}=p^*$.  As $p$ is $A$-large, $R$ is uncountable.  Now, for each finite $\mbar$ from $M$,
there are only countably many complete $q\in S_{at}(\mbar)$, and if some $q\in S_{at}(\mbar)$ is $A$-small, then only countably many $\pbar\in R$ extend $q$.
As $M$ is countable, there are only countably many $\mbar$, hence all but countably many $\pbar\in R$ satisfy $\pbar\mr{\mbar}$ $A$-large for every $\mbar$.
Further, again since $M$ is countable, at most countably many $\pbar\in R$ have restrictions to $A$ that are realized in $M$.  Thus, all but countably many $\pbar\in R$ are 
$A$-perfect.
\endproof

\subsection{A transfer result}
In this brief subsection we state a transfer result that follows immediately by Keisler's completeness theorem for the logic  $L_{\omega_1,\omega}(Q)$, given in \cite{Keisler70}.
Recall that $L_{\omega_1,\omega}(Q)$ is the logic obtained by taking the (usual) set of atomic $L$ formulas and closing under boolean combinations, existential quantification,
the `$Q$-quantifier,' i.e., if $\theta(y,\xbar)$ is a formula, then so is $Qy\theta(y,\xbar)$; and countable conjunctions of formulas involving a finite set of free variables,
i.e., if $\{\psi_i(\xbar):i\in\omega\}$ is a set of formulas, then so is $\bigwedge_{i\in\omega} \psi_i(\xbar)$.  We are only interested in {\em standard interpretations} of these formulas,
i.e., $M\models \bigwedge_{i\in\omega}\psi_i(\abar)$ if and only if $M\models\psi_i(\abar)$ for every $i\in\omega$; and $M\models Qy\theta(y,\abar)$ if and only if the solution
set $\theta(M,\abar)$ is uncountable.

Throughout the discussion let $ZFC^*$ denote a sufficiently large, finite subset of the ZFC axioms.  In the notation of \cite{URL},  Proposition~\ref{trans} states that sentences of $L_{\omega_1,\omega}(Q)$ are {\em grounded.}

\begin{Proposition}  \label{trans}  Suppose $L$ is a countable language, and $\Phi\in L_{\omega_1,\omega}(Q)$ are given.
There is a sufficiently large, finite subset $ZFC^*$ of $ZFC$ such that
IF there is a countable, transitive model $(\B,\epsilon)\models ZFC^*$ with $L,\Phi\in \B$ and
$$(\B,\epsilon)\models \ \hbox{`There is $M\models\Phi$ and $|M|=\aleph_1$'}$$
THEN (in $V$!) there is $N\models\Phi$ and $|N|=\aleph_1$.
\end{Proposition}


\bp  This follows immediately from Keiser's completeness theorem for $L_{\omega_1,\omega}$, given that provability is absolute between transitive models of set theory.
More modern, `constructive' proofs can be found in \cite{BLarson} and \cite{BLS}.  These use the existence $\B$-normal ultrafilters.  Given an arbitrary language $L^*\in\B$
and any countable $L^*$-structure $(\B,E,\dots)$ where  the reduct $(\B,E)$ is an $\omega$-model of $ZFC^*$, for any $\B$-normal ultrafilter $\U$, the ultrapower
$Ult(\B,\U)$ is a countable, $\omega$-model that is an $L^*$-elementary extension of $(\B,E,\dots)$.  It has the additional property that for any $L^*$-definable subset
$D$, $D^{Ult(\B,\U)}$ properly extends $D^{\B}$ if and only if $(\B,E,\dots)\models `D$ is uncountable'.

Using this, one constructs (in $V$!) a continuous, $L^*$-elementary $\omega_1$-sequence  $\<\B_\alpha:\alpha<\omega_1\>$ of $\omega$-models, where each 
$\B_{\alpha+1}=Ult(\B_\alpha,\U_\alpha)$.  Then the interpretation $M^{\CC}$ where $\CC=\bigcup_{\alpha\in\omega_1} \B_\alpha$ will be a suitable choice of $N$.
More details of this construction are given in \cite{BLarson} or \cite{BLS}.

\subsection{A configuration arising from non-$\pcl$-smallness}

The goal of this subsection is to prove the following Proposition, the data from which will be used throughout Section~\ref{three}.

\begin{Proposition}  \label{config} Assume $T$ is a countable, complete theory for which $\At_T$ has an uncountable atomic model, but is not $\pcl$-small.  
Then there are a countable, atomic  $M^*\in\At_T$, finite sequences $\abar^*\subseteq\bbar^*\subseteq M^*$, and complete 1-types $\{r_j(x,\bbar^*):j\in\omega\}$
such that, letting $D^*=\pcl_{M^*}(\abar^*)$, $A_n=\bigcup\{r_j(M^*,\bbar^*):j<n\}$ and $A^*=\bigcup\{A_n:n\in\omega\}$ we have:
\begin{enumerate}
\item  $A^*\subseteq D^*$;
\item  $S^+_{at}(A_n,M^*)$ is countable for every $n\in\omega$; but
\item  $S^+_{at}(A^*,M^*)$ is uncountable.
\end{enumerate}
\end{Proposition}

\bp  Fix any countable, atomic $M^*\in\At_T$.  Using Proposition~\ref{relate} and the non-$\pcl$-smallness of $\At_T$,  choose a finite tuple $\abar^*\subseteq M^*$
such that $S^+_{at}(D^*,M^*)$ is uncountable, where
 $D^*=\pcl_{M^*}(\abar^*)\subseteq M^*$.

 Fix any finite tuple $\bbar\supseteq\abar^*$ from $M^*$ and look at the complete 1-types $\QQ_{\bbar}:=\{r\in S_{at}(\bbar)$ such that $r(M^*)\subseteq D^*\}$.
 These types visibly induce a  partition $D^*$, and it is easily seen that if $\bbar'\supseteq\bbar$, the partition induced by $\bbar'$ refines the partition induced by $\bbar$.
 Let $\QQ:=\bigcup\{\QQ_{\bbar}:\abar^*\subseteq\bbar\subseteq M^*\}$.
 

%
%
%
%
 
Define a rank function $\rk:\QQ\rightarrow ON\cup\{\infty\}$  as follows:
 \begin{itemize}
 \item  $\rk(c/\bbar)\ge 0$ if and only if $\tp(c/\bbar)\in\QQ$;
 \item  $\rk(c/\bbar)\ge 1$ if and only if $\tp(c/\bbar)\in\QQ$ and there are infinitely many $c'\in D^*$ realizing $\tp(c/D^*)$; and
 \item  for an ordinal $\alpha\ge 2$, $\rk(c/\bbar)\ge \alpha$ if and only if for every $\beta<\alpha$ and every $\bbar'$ from $M^*$,
 there is $c'\in D^*$ realizing $\tp(c/\bbar)$ such that $\rk(c'/\bbar\bbar')\ge\beta$.
 \item  $\rk(c/\bbar)=\alpha$ if and only if $\rk(c/\bbar)\ge\alpha$ but $\rk(c/\bbar)\not\ge\alpha+1$.
 \end{itemize}
 
 \medskip\noindent{\bf Claim 1.}  For every $r\in\QQ$, $\rk(r)$ is a countable ordinal.
 
 \medskip
 
 \bp  Assume by way of contradiction that $\rk(c/\bbar)\ge\omega_1$ for some type $c/\bbar$.  Then, for any $\bbar'$ from $M$, as $D^*$ is countable, there is
 an element $c'\in D^*$ such that $\rk(c'/\bbar\bbar')\ge\beta$ for uncountably many $\beta$'s, hence $\rk(c'/\bbar\bbar')\ge\omega_1$ as well.
 Using this idea, if we let $\<\bbar_n:n\in\omega\>$ be an increasing sequence of finite sequences from $M^*$ whose union is all of $M^*$, 
 then we can find a sequence $\<c_n:n\in\omega\>$ of elements from $D^*$ such that, for each $n$, $\rk(c_n/\bbar_n)\ge\omega_1$ 
 and $\tp(c_n/\bbar_n)\subseteq\tp(c_{n+1}/\bbar_{n+1})$.  
The union of these 1-types yields a complete, atomic 1-type $q\in S_{at}(M^*)$ all of whose realizations are in $\pcl_{M^*}(\abar)$.  
However, since the type asserting that `$x=c$' has rank 0 for each $c\in D^*$,  $q$ is omitted in $M^*$.  To obtain a contradiction, choose a realization $e$ of $q$ and,
as $M^*e$ is a countable, atomic set, construct a countable, elementary extension $M'\succeq M^*$ with $e\in M'$.  But now, $q$ implies that $e\in\pcl_{M'}(\abar)$,
yet this is contradicted by the fact that $M^*$ contains $\abar$ but not $e$.
\endproof

As notation, for a subset $\SS\subseteq\QQ_{\bbar}$, let $A_{\SS}=\bigcup\{r(M^*):r\in\SS\}$, which is always a subset of $D^*$.
Define the set of `candidates' as
$$\CC=\{(\SS,\bbar):\bbar\supseteq\abar^*, \SS\subseteq\QQ_{\bbar}, \ \hbox{and}\ S^+_{at}(A_{\SS},M^*)\ \hbox{uncountable}\}$$
Note that $\CC$ is non-empty as $(\SS_0,\abar^*)\in\CC$, where $\SS_0$ is an enumeration of all the complete, pseudo-algebraic types over $\abar^*$.
Among all candidates, choose $(\SS^*,\bbar^*)\in\CC$ such that
$$\alpha^*:=\sup\{rk(r)+1:r\in\SS^*\}$$ is as small as possible.  
Enumerate $\SS^*=\{r_j:j\in\omega\}$ and put $A^*:=A_{\SS^*}$ and $A_n:=\bigcup\{r_j(M^*,\bbar^*):j<n\}$ for each $n\in\omega$.
As Clauses~(1) and (3) are immediate, it suffices to prove the following Claim:

\medskip\noindent{\bf Claim 2.}  For each $n\in\omega$, $S^+_{at}(A_n,M^*)$ is countable.
 
 \medskip
 
 \bp Fix any $n\in\omega$.  First, note that if $\rk(r_j)=0$ for every $j<n$, then $A_n$ would be finite, which would imply $S_{at}(A_n)$ is countable.  As 
 $S_{at}(A_n)$ contains $S^+_{at}(A_n,M^*)$, the result follows.
 
 Now assume $rk(r_j)>0$ for at least one $j<n$.  Let $\beta:=\max\{rk(r_j):j<n\}$ and let $F=\{j<n:rk(r_j)=\beta\}$.
 Clearly, $\beta<\alpha^*$.
For each $j\in F$, as $\beta>0$ but $rk(r_j)\not\ge \beta+1$, there is a finite tuple $\bbar_j$ such that
$\rk(c/\bbar^*\bbar_j)<\beta$ for all $c\in r_j(M^*)$.  

Let $\bbar'$ be the concatenation of $\bbar^*$ with each $\bbar_j$ for $j\in F$ and let
$$\SS':=\{r'\in\QQ_{\bbar'}:r'\ \hbox{extends some $r_j$ with $j<n$}\}$$

 \medskip\noindent{\bf Subclaim.}  $\rk(r')<\beta$ for every $r'\in\SS'$.

 \medskip
 
 \bp Fix $r'\in \SS'$ and choose $c\in r'(M^*,\bbar')$.  There are two cases.  On one hand, if $r'$ extends some $r_j$ with $j\in F$, then
 $\rk(c/\bbar')\le\rk(c/\bbar^*\bbar_j)<\beta$.  On the other hand, if $r'$ extends some $r_j$ with $r_j\not\in F$, then as $\rk(r_j)<\beta$,
 $\rk(c/\bbar')\le\rk(c/\bbar^*)<\beta$.
 \endproof

Clearly $A_{\SS'}=A_n$, so $S^+_{at}(A_n,M^*)=S^+_{at}(A_{\SS'},M^*)$.  
Thus, if $S^+_{at}(A_n,M^*)$ were uncountable, then $(\SS',\bbar')$ would be a candidate, i.e., an element of $\CC$.
But, as $\beta<\alpha^*$, this is impossible by the Subclaim and the minimality of $\alpha^*$.
\endproof

  \section{A family of $2^{\aleph_0}$ atomic models of size $\aleph_1$}  \label{three}
 
 Throughout the whole of this section, we assume that $T$ is a complete theory in a countable language for which $\At_T$ has an uncountable atomic model, but is
 not $\pcl$-small.  Appealing to Proposition~\ref{config},  
 \begin{quotation}
 
 \noindent {\bf Fix, for the whole of this section, a countable atomic model ${\bf M^*}$, tuples  ${\bf \abar^*\subseteq\bbar^*\subseteq M^*}$ and sets $A^*$ and $A_n$ for each
 $n\in\omega$ as in Proposition~\ref{config}.}
  \end{quotation}
 
%
%

We work with this fixed configuration for the whole of this section and, in Subsection~\ref{massp} eventually prove:
 
 \begin{Proposition}  \label{many}  There is a family $\{(N_\eta,\bbar^*):\eta\in 2^\omega\}$ of atomic models of $T$, each of size $\aleph_1$,
 that are pairwise non-isomorphic over $\bbar^*$.
 \end{Proposition}

 
 \subsection{Colorings of models realizing many types over $A^*$}

\begin{Definition} {\em  Call a structure $(N,\bbar^*)$ {\em rich} if $N\in \At_T$ has size $\aleph_1$,  $M^*\preceq N$, and $N$ realizes uncountably many 
1-types over $A^*$.
}
\end{Definition}

 
 \begin{Lemma}  \label{count}
For each $n\in\omega$, a rich $(N,\bbar^*)$ realizes only countably many distinct 1-types over $A_n$.
 \end{Lemma}
 
 \bp  Fix any $(N,\bbar^*)$ and $n<\omega$ as above.  If $\{c_i:i\in\omega_1\}$ realize distinct types over $A_n$, then the types $\{\tp_N(c_i/M^*):i\in\omega_1\}$
 would be distinct, contradicting $S^+_{at}(A_n,M^*)$ countable.
 \endproof
%
%
%
 
 How can we tell whether rich  structures are non-isomorphic?    
We introduce the notion of $\U$-colorings and Corollary~\ref{criterion} gives a sufficient condition.

\begin{Definition} \label{color} {\em Fix a subset $\U\subseteq\omega$ and a rich $(N,\bbar^*)$.
\begin{itemize}  
\item  For elements $d,d'\in N$, define the {\em splitting number} $\spl(d,d')\in(\omega+1)$
to be the least $k<\omega$ such that $\tp(d/A_k)\neq\tp(d'/A_k)$ if such exists; and $\spl(d,d')=\omega$ if 
$\tp(d/A^*)=\tp(d'/A^*)$.

\item  A {\em $\U$-coloring of  a rich $(N,\bbar^*)$} is a function
$$c:N\rightarrow \omega$$
such that for all pairs $d,d'\in N$, at least one of the following hold:
\begin{enumerate}
\item   $\tp(d/A^*)=\tp(d'/A^*)$; or
\item $c(d)\neq c(d')$; or
\item  $\spl(d,d')\in\U$.
\end{enumerate}
\item  The {\em color filter} $\F(N,\bbar^*):=\{\U\subseteq\omega:$ a $\U$-coloring of $(N,\bbar^*)$ exists$\}$.
\end{itemize}
}
\end{Definition}

\begin{Lemma}  \label{test}   Fix a rich $(N,\bbar^*)$.  Then:
\begin{enumerate}
\item  $\F(N,\bbar^*)$ is a filter; 
\item  $\F(N,\bbar^*)$ contains the cofinite subsets of $\omega$; but
\item  No finite $\U\subseteq\omega$ is in $\F(N,\bbar^*)$. 
\end{enumerate}
\end{Lemma}

\bp  (1)  First, note that if $\U\subseteq\U'\subseteq\omega$, then every $\U$-coloring $c$ is also a $\U'$-coloring.  Thus, $\F(N,\bbar^*)$ is upward closed.
Next, suppose $\U_1\in\F(N,\bbar^*)$ via the coloring $c_1:N\rightarrow\omega$ and $\U_2\in\F(N,\abar^*\bbar^*)$ via the coloring $c_2:N\rightarrow\omega$.
Fix any bijection $t:\omega\times\omega\rightarrow\omega$.  It is easily checked that $c^*:N\rightarrow\omega$ defined by $c^*(d)=t(c_1(d),c_2(d))$ is a 
$\U_1\cap \U_2$-coloring of $(N,\bbar^*)$.  Thus, $\U_1\cap\U_2\in\F(N,\bbar^*)$.  So $\F(N,\bbar^*)$ is a filter.

(2)   As $\F(N,\bbar^*)$ is a filter, it  suffices to  show 
$(\omega-n)\in\F(N,\bbar^*)$ for each $n\in\omega$.  So fix such an $n$.  By Lemma~\ref{count}, $N$ realizes at most countably many types over $A_n$.
Thus, we can  produce a map $c:N\rightarrow\omega$ such that $c(d)=c(d')$ if and only if $\tp(d/A_n)=\tp(d'/A_n)$.  As any such $c$ is an $(\omega-n)$-coloring,
$(\omega-n)\in\F(N,\bbar^*)$.

(3)  It suffices to show that no $n=\{0,\dots,n-1\}$ is in $\F(N,\bbar^*)$.  To see this, let $c:N\rightarrow\omega$ be an arbitrary map.  
We will show that $c$ is not an $\{0,\dots,n-1\}$-coloring.  As 
$N$ realizes $\aleph_1$ distinct types over $A^*$, there is some $m^*\in\omega$ and an uncountable subset 
$\{d_\alpha:\alpha<\omega_1\}\subseteq N$ that realize distinct types over $A^*$, yet $c(d_\alpha)=m^*$ for each $\alpha$.
However, as $N$ realizes only countably many types over $A_n$, there are $\alpha\neq\beta$ such that $n\le \spl(d_\alpha,d_\beta)<\omega$.
Thus, $c$ is not an $\{0,\dots,n-1\}$-coloring.
\endproof

We close with a sufficient condition for non-isomorphism of rich models.

\begin{Corollary}  \label{criterion} Suppose that for $\ell=1,2$, $(N_\ell,\bbar^*)$ is a $\U_\ell$-colored rich model, and $\U_1\cap\U_2$ is finite.  Then there is
no isomorphism $f:N_1\rightarrow N_2$ fixing $\bbar^*$ pointwise.
\end{Corollary}

\bp  If there were such an isomorphism, then $(N_2,\bbar^*)$ would be both $\U_1$-colored and $\U_2$-colored.  Thus, both $\U_1,\U_2\in\F(N_2,\bbar^*)$,
which contradicts Lemma~\ref{test}.
\endproof

\subsection{Constructing a colored rich model via forcing}

Arguing as in the proof of Proposition~\ref{relate}, from the data of Lemma~\ref{config} we can construct a rich $(N,\bbar^*)$
as the union of a continuous, elementary chain $\<M_\alpha:\alpha\in \omega_1\>$
of countable, atomic models with $M_0=M^*$ such that, for each $\alpha\in\omega_1$ there is a distinguished $b_\alpha\in M_{\alpha+1}$ such that $\tp(b_\alpha/A^*)$ is
omitted in $M_\alpha$.    

Our goal is to construct a sufficiently generic rich $(N,\bbar^*)$, along with a coloring $c:N\rightarrow(\omega+1)$ via forcing.
Our forcing $(\Q,\le_{\Q})$ encodes finite approximations of such an $(N,\bbar^*)$ and $c$.  A fundamental building block is the notion of
a {\em striated type} over a finite subset $\abar$  satisfying $\bbar^*\subseteq\abar\subseteq M^*$.  
 As an atomic type over a finite subset is generated by a complete formula,
we use the terms interchangeably.

\begin{Definition}  \label{striated}  {\em  Choose a finite tuple $\abar$ with $\bbar^*\subseteq\abar\subseteq M^*$.  A {\em striated type over $\abar$} is a complete formula
$\theta(\xbar)\in S_{at}(\abar)$ whose variables are partitioned as $\xbar=\<\xbar_j:j<\ell\>$ where, for each $j$,   $\xbar_j=\<x_{j,n}:n< n(j)\>$ is an $n(j)$-tuple of variable symbols
that satisfy $\tp(x_{j,0}/\abar\cup\{\xbar_i:i<j\})$ is $A^*$-large.  The integer $\ell$ is the {\em length} of the striated type.

A {\em simple realization} of a striated type $\theta(\xbar)$  of length $\ell$ is a sequence  $\bbar=\<\bbar_j:j<\ell\>$ of tuples from $M^*$ such that $M^*\models \theta(\bbar)$.
A {\em perfect chain realization} of  $\theta(\xbar)$ is a pair $(\Mbar,\bbar)$, consisting of  a chain $M_0\preceq M_1\preceq M_{\ell-1}\preceq M^*$
of  $\ell$ elementary submodels  of $M^*$ and a simple realization $\bbar=\<\bbar_j:j<\ell\>$  from $M^*$ that satisfy:  For each $j<\ell$,
\begin{enumerate}
\item  $\abar\cup\{\bbar_i:i<j\}\subseteq M_j$; and
\item  $\tp(b_{j,0}/M_j)$ is $A^*$-perfect (see Definition~\ref{perfectx}).
\end{enumerate}
}
\end{Definition}

\begin{Lemma}  \label{existsg}  Every striated type $\theta(\xbar)\in S_{at}(\abar)$ has a perfect chain realization.
\end{Lemma}

\bp  We argue by induction on $\ell$, the length of the striation.
For striations of length zero there is nothing to prove, so assume the Lemma holds for striated types of length $\ell$
and choose an $(\ell+1)$-striation $\theta(\xbar)\in S_{at}(\abar)$.  Let $\theta\mr{\ell}$ be the truncation of $\theta$ to the variables
$\xbar\mr{\ell}=\<\xbar_j:j<\ell\>$.  As $\theta\mr{\ell}$ is clearly an $\ell$-striation,  it has a perfect chain realization, i.e.,  a chain $M_0\preceq M_1\preceq M_{\ell-1}\preceq M^*$
and a tuple $\bbar=\<\bbar_j:j<\ell\>$ from $M^*$  realizing $\theta\mr{\ell}$ such that $\abar\cup\{\bbar_i:i<j\}\subseteq M_j$ and $\tp(b_{j,0}/M_j)$ is $A^*$-perfect for each
$j<\ell$.

Now, since $\tp(x_{\ell,0}/\abar\bbar)$ is $A^*$-large, by applying Proposition~\ref{perfect} there is  an $A^*$-perfect type $\pbar\in S_{at}(M^*)$
 (in a single variable $x_{\ell,0}$) extending $\tp(x_{\ell,0}/\abar\bbar)$.
 Choose a countable, atomic $N\succeq M^*$ and $e\in N$ realizing $\pbar$.  As $N$ and $M^*$ are both countable and atomic,
choose an isomorphism $f:N\rightarrow M^*$ that fixes $\abar\bbar$ pointwise.  
Then $f(M_0)\preceq f(M_1)\preceq \dots f(M_{\ell-1})\preceq f(M^*)\preceq M^*$ is a chain.  Let $b_{\ell,0}:=f(e)$ and choose $\<b_{\ell,1}\dots,b_{\ell,n(\ell)-1}\>$ arbitrarily from
$M^*$ so that, letting $\bbar_\ell=\<\bbar_{\ell,n}:n<n(\ell)\>$, $\bbar\frown \bbar_\ell$ realizes $\theta(\xbar)$.  This chain and this sequence form a perfect chain realization of $\theta$.
\endproof

The following Lemma is immediate, and indicates the advantage of working with $A^*$-perfect types.

\begin{Lemma}  \label{extend} Let $(\Mbar,\bbar)$ be any perfect chain realization of a striated type 
$\theta(\xbar)\in S_{at}(\abar)$.  Then for every $\cbar\subseteq M_0$, $\tp(\bbar/\abar\cbar)\in S_{at}(\abar\cbar)$ is a striated type extending $\theta(\xbar)$,
and $(\Mbar,\bbar)$ is a perfect chain realization of it.
\end{Lemma}

The Lemma below, whose proof simply amounts to unpacking definitions, demonstrate that striated types are rather malleable.

\begin{Lemma}  \label{close}  
\begin{enumerate}
\item  If $\tp(\cbar/\abar)$ is a striated type of length $k$ and $\tp(\dbar/\abar\cbar)$ is a striated type of length $\ell$, then
$\tp(\cbar\dbar/\abar)$ is a striated type of length $k+\ell$.
\item  Suppose $\tp(\bbar/\abar)$ is a striated type of length $\ell$ and $k<\ell$.  Let $\bbar_{<k}$ and $\bbar_{\ge k}$ be the induced partition of $\bbar$.
Then $\tp(\bbar_{<k}/\abar)$ is a striated type of length $\ell$ and $\tp(\bbar_{\ge k}/\abar\bbar_{<k})$ is a striated type of length $(\ell-k)$.  
Moreover, if $(\Mbar,\bbar)$ is a perfect chain realization of $\tp(\bbar/\abar)$, then $(\Mbar_{<k},\bbar_{<k})$ is a perfect chain realization of
$\tp(\bbar_{<k}/\abar)$ and $(\Mbar_{\ge k},\bbar_{\ge k})$ is a perfect chain realization of $\tp(\bbar_{\ge k}/\abar\bbar_{<k})$.
\end{enumerate}
\end{Lemma}

We begin by defining a partial order  $(\Q_0,\le_{\Q_0})$ of `preconditions'.  Then our forcing $(\Q,\le_{\Q})$
will  be a dense suborder of these preconditions. 

\begin{Definition}  \label{forcing}  {\em  $\Q_0$ is the set of all $\pp=(\aabar_{\pp}, u_{\pp}, \nbar_{\pp}, \theta_{\pp}(\xbar_{\pp}),k_{\pp},\U_{\pp},c_{\pp})$, where
\begin{enumerate}
\item  $\aabar_{\pp}$ is a finite subset of $M^*$ containing $\bbar^*$;
\item $u_{\pp}$ is a finite subset of $\omega_1$;
\item  $\nbar_{\pp}=\<n_t:t\in u_{\pp}\>$ is a sequence of positive integers;
\item  $\xbar_{\pp}=\<\xbar_{t,\pp}:t\in u_{\pp}\>$, where each $\xbar_{t,\pp}=\<x_{t,n}: n< \nbar_{t}\>$ is a finite sequence from the  set $X=\{x_{t,n}:t\in\omega_1,n\in\omega\}$ of variable symbols;
\item  $\theta_{\pp}(\xbar_{\pp})\in S_{at}(\aabar_{\pp})$ is a striated type of length $|u_{\pp}|$ (see Definition~\ref{striated});
\item  $k_{\pp}\in\omega$;
\item  $\U_{\pp}\subseteq k_{\pp}=\{0,\dots,k_{\pp}-1\}$;
\item  $c_{\pp}:\xbar_{\pp}\rightarrow\omega$ is a function
such that 
for all pairs $x_{t,n},x_{s,m}$ from $\xbar_{\pp}$ with $c_{\pp}(x_{t,n})=c_{\pp}(x_{s,m})$
\begin{enumerate}  \item {\bf either} 
$\spl(b_{t,n},b_{s,m})\ge k_{\pp}$ for all perfect chain realizations $(\Mbar,\bbar)$ of $\theta_{\pp}(\xbar_{\pp})$;
\item {\bf or}  there is some $k\in\U_{\pp}$ such that $\spl(b_{t,n},b_{s,m})=k$ for all perfect chain realizations $(\Mbar,\bbar)$ of $\theta_{\pp}(\xbar_{\pp})$.
\end{enumerate}

\end{enumerate}
}
\end{Definition}

We order elements of $\Q_0$ by: $\pp\le_{\Q_0}\qq$ if and only if
\begin{itemize}
\item  $\aabar_{\pp}\subseteq \aabar_{\qq}$;
\item  $u_{\pp}\subseteq u_{\qq}$ and $n_{t,\pp}\le n_{t,\qq}$ for all $t\in u_{\pp}$, hence $\xbar_{\pp}$ is a subsequence of $\xbar_{\qq}$;
\item  $\theta_{\qq}(\xbar_{\qq})\vdash \theta_{\pp}(\xbar_{\pp})$;
\item  $k_{\pp}\le k_{\qq}$;
\item  $\U_{\pp}=\U_{\qq}\cap k_{\pp}$ (hence, for $j< k_{\pp}$, $j\in\U_{\pp}$ if and only if $j\in \U_{\qq}$);
\item  $c_{\pp}=c_{\qq}\mr{\xbar_{\pp}}$.
\end{itemize}

Visibly, $(\Q_0,\le_{\Q_0})$ is a partial order.  
Call a precondition $\pp\in\Q_0$ {\em  unarily decided} if, for every
$x_{t,n}\in \xbar_{\pp}$, $p(\xbar_{\pp})$ determines a type in $S^+_{at}(A_{k_{\pp}},M^*)$ (see Definition~\ref{lies}).    That the  unarily decided 
preconditions are dense follows easily from the fact that
$S^+_{at}(A_{k_{\pp}},M^*)$ is countable.

\begin{Lemma}  The set $\{\pp\in\Q_0:\pp$ is  unarily decided$\}$ is dense in $(\Q_0,\le_{\Q_0})$.  Moreover, given any $\pp\in \Q_0$, there is a unarily decided
$\qq\ge_{\Q_0}\pp$ with  $\xbar_{\qq}=\xbar_{\pp}$  and $k_{\qq}=k_{\pp}$ (hence $\U_{\qq}=\U_{\pp}$).
\end{Lemma}

\bp  Fix $\pp\in\Q_0$ and  let $k:=k_{\pp}$.
Arguing by induction on the size of the finite set $\xbar_{\pp}$, it is enough to strengthen $p(x_{t,n})$ individually for each $x_{t,n}\in \xbar_{\pp}$.
So fix $x_{t,n}\in \xbar_{\pp}$.  By Corollary~\ref{countable}  there is an $\abar'\supseteq \abar_{\pp}$ and a 1-type $q_1(x_{t,n})\in S_{at}(\abar')$
extending $\tp(x_{t,n}/\abar_{\pp})$ that determines a type in $S_{at}^+(A_{k_{\pp}},M^*)$.  Then, using Lemma~\ref{close}(1) we can choose a striated type
$p'(\xbar_{\pp})\in S_{at}(\abar')$ extending $p(\xbar_{\pp})\cup q_1$.

We iterate the above procedure for each of the (finitely many) elements of $\xbar_{\pp}$.  We then get a unarily decided precondition
 $\pp'\ge_{\Q_0} \pp$  whose type $p'(\xbar_{\pp})$ still has the same free variables, and each of $k_{\pp}$, $\U_{\pp}$, $c_{\pp}$  are unchanged.
\endproof

Next, call a precondition $\pp\in\Q_0$ {\em fully decided} if, it is unarily decided and, for each pair $x_{t,n},x_{s,m}$ from $\xbar_{\pp}$ with $c_{\pp}(x_{t,n})=c_{\pp}(x_{s,m})$, if   
$\spl(b_{t,n},b_{s,m})\ge k_{\pp}$
for some perfect chain realization $(\Mbar,\bbar)$, then $\tp(b_{t,n}/A^*)=\tp(b_{s,m}/A^*)$ for all perfect chain realizations $(\Mbar,\bbar)$ of $\theta_{\pp}(\xbar_{\pp})$.

\begin{Lemma} \label{full} The set $\{\pp\in\Q_0:\pp$ is  fully decided$\}$ is dense in $(\Q_0,\le_{\Q_0})$.  Moreover, given any $\pp\in \Q_0$, there is a fully decided
$\qq\ge_{\Q_0}\pp$ with  $\xbar_{\qq}=\xbar_{\pp}$.
\end{Lemma}

\bp  It suffices to  handle each pair $x_{t,n},x_{s,m}$ from $\xbar_{\pp}$ with $c(x_{t,n})=c(x_{s,m})$ separately.  Given such a pair, suppose there is some
perfect chain realization $(\Mbar,\bbar)$ of $\theta(\xbar_{\pp})\in S_{at}(\aabar_{\pp})$ with $k_{\pp}\le \spl(b_{t,n},b_{s,m})<\omega$.  Among all such perfect chain realizations,
choose one that minimizes $k^*=\spl(b_{t,n},b_{s,m})$.  
Choose a formula $\phi(x,\cbar)$ with $\cbar$ from $A_{k^*+1}$ witnessing that $\tp(b_{t,n}/A_{k^*+1})\neq\tp(b_{s,m}/A_{k^*+1})$.  
As $A_{k^*+1}\subseteq M_0$, by applying Lemma~\ref{extend}, let $\theta^*(\xbar_{\pp})$ be a complete formula over $\aabar_{\pp}\cbar$ isolating $\tp(\bbar/\aabar_{\pp}\cbar)$.
Form the precondition $\pp'\in\Q_0$ by  putting $\aabar_{\pp'}=\aabar_{\pp}\cbar$; 
$\theta_{\pp'}=\theta^*$; $k_{\pp'}=k^*+1$; and $\U_{\pp'}=\U_{\pp}\cup\{k^*\}$; while leaving
$\xbar_{\pp}$ and $c_{\pp}$ unchanged.
It is evident that $\spl(b'_{t,n},b'_{s,m})=k^*\in\U_{\pp'}$ for all perfect chain realizations $(\Mbar,\bbar')$ of $\theta_{\pp'}$.  Continuing this process for each of the (finitely many)
relevant pairs gives us a fully decided extension of $\pp$.
\endproof

\begin{Definition}  The forcing $(\Q,\le_\Q)$ is the set of fully decided $\pp\in\Q_0$ with the inherited order.
\end{Definition}

\begin{Lemma}  The forcing $(\Q,\le_{\Q})$ has the  countable chain condition (c.c.c.).
\end{Lemma}

\bp  Suppose $\{\pp_i:i\in\omega_1\}$ is an uncountable subset of $\Q$.  In light of Lemma~\ref{full}, it suffices to find $i\neq j$ for which there is some precondition
$\qq\in\Q_0$ satisfying $\pp_i\le_{\Q_0} \qq$ and $\pp_j\le_{\Q_0} \qq$.  First, by the $\Delta$-system lemma  applied to the finite
sets $\{u_{\pp_i}\}$, we may assume that $|u_{\pp_i}|$ is constant and there is some fixed $u^*$ that is an initial segment of each $u_{\pp_i}$ and, moreover, whenever $i<j$,
every element of $(u_{\pp_i}\setminus u^*)$ is less than every element of $(u_{\pp_j}\setminus u^*)$.
By further trimming, but preserving uncountability, we may  assume that
the integer $k_{\pp}$, the subset $\U_{\pp}\subseteq k_{\pp}$, and the parameter $\aabar_{\pp}$ remain constant.  
As notation, for $i<j$, let $f:u_{\pp_i}\rightarrow u_{\pp_j}$ be the unique order-preserving bijection.   We may additionally assume that $n_{\pp_i}(t)=n_{\pp_j}(f(t))$, hence
$f$ has a natural extension (also called $f$)$:\xbar_{\pp_i}\rightarrow\xbar_{\pp_j}$ given by $f(x_{t,n})=x_{f(t),n}$.  With this identification, we may assume 
$\theta_{\pp_i}(\xbar_{\pp_i})=\theta_{\pp_j}(f(\xbar_{\pp_i}))$.  As well, we may
also assume
$\tp(x_{t,n}/A_{k_{\pp}})=\tp(x_{f(t),n}/A_{k_{\pp}})$ for every $x_{t,n}\in\xbar_{\pp_i}$.  As well, the colorings match up as well, i.e., $c(x_{t,n})=x_{f(t),n}$.

Now fix $i<j$.  Define $\qq$ by $k_{\qq}:=k_{\pp}$; $\U_{\qq}:=\U_{\pp}$; and $\aabar_\qq:=\aabar_{\pp}$ (the common values).  Let $u_{\qq}:=u_{\pp_i}\cup u_{\pp_j}$, 
and, for $t\in u_{\pp_i}$, $n_{t,\qq}=n_{t,\pp_i}$ while $n_{t,\qq}=n_{t,\pp_j}$ for $t\in u_{\pp_j}$.  To produce the striated type $\theta_{\qq}\in S_{at}(\aabar_{\qq})$, first choose
a perfect chain realization $(\Mbar,\bbar)$ of $\theta_{\pp_i}(\xbar_{\pp_i})$.  Say  $|u_{\pp_i}|=\ell=|u_{\pp_j}|$, while $|u^*|=k<\ell$.  
By Lemma~\ref{close}(2), $\tp(\bbar_{<k}/\aabar_{\pp})$ is a striated type of length $k$ and 
$(\Mbar_{\ge k},\bbar_{\ge k})$ is a perfect chain realization of the striated type $\tp(\bbar_{\ge k}/\aabar_{\pp}\bbar_{<k})$ of length $(\ell-k)$.
Choose $\dbar$ from $M_k$ such that $\tp(\dbar/\aabar_{\pp}\bbar_{<k})=\tp(\bbar_{\ge k}/\aabar_{\pp}\bbar_{<k})$.
Then by Lemma~\ref{extend} (with $M_k$ playing the role of $M_0$ there), $(\Mbar_{\ge k},\bbar_{\ge k})$ is a perfect chain realization of the striated type
$\tp(\bbar_{\ge k}/\aabar_{\pp}\bbar_{<k}\dbar)$.  So, by Lemma~\ref{close}(1), $\tp(\dbar\bbar_{\ge k}/\aabar_{\pp}\bbar_{<k})$ is a striated type of length $2(\ell-k)$.
Thus, a second application of Lemma~\ref{close}(1) implies that $\tp(\bbar_{<k}\dbar\bbar_{\ge k}/\aabar_{\pp})$ is a striated type of length $2\ell-k$.
Let $\theta_{\qq}$ be a complete formula over $\aabar_{\pp}$ generating this type.  


In order to show that $\qq$ is a precondition (i.e., an element of $\Q_0$) only Clause~(8) requires an argument.
Fix any  $x_{t,n},x_{s,m}$ in $\xbar_{\qq}$ with $c_{\qq}(x_{t,n})=c_{\qq}(x_{s,m})$.  As both $\pp_i,\pp_j\in\Q_0$, the verification is immediate if $\{t,s\}$ is a subset of either
$u_{\pp_i}$ or $u_{\pp_j}$, so assume otherwise.  By symmetry, assume $t\in u_{\pp_i}-u^*$ and $s\in u_{\pp_j}-u^*$.  The point is that by our trimming,
$x_{f(t),n}\in \xbar_{\pp_j}$, $c_{\pp_j}(x_{f(t),n})=c_{\pp_i}(x_{t,n})$, and $\tp(x_{t,n}/A_{k_{\pp}})=\tp(x_{f(t),n}/A_{k_{\pp}})$.  
There are now two cases:  First, if $\tp(x_{f(t),n}/A^*)=\tp(x_{s,m}/A^*)$, then it follows that $\tp(x_{t,n}/A_{k_{\pp}})=\tp(x_{s,m}/A_{k_{\pp}})$, hence
$\spl(e_{t,n},e_{s,m})\ge k_{\pp}$ for any perfect chain realization $(\Nbar,\ebar)$ of $\theta_{\qq}$.
On the other hand, if $\theta_{\pp_j}$ `says' $\spl(x_{f(t),n},x_{s,m})=k\in\U_{\pp}$, then $\theta_{\qq}$ `says' $\spl(x_{t,n},x_{s,m})=k\in\U_{\qq}$ as well.
Thus, $\qq\in\Q_0$, which suffices by Lemma~\ref{full}.
\endproof

\begin{Lemma}  Each of the following sets are dense and open in $(\Q,\le_{\Q})$.
\begin{enumerate}
\item  For every $t\in\omega_1$, $D_t=\{\pp\in\Q: t\in u_{\pp}\}$;
\item  For every $(t,n)\in\omega_1\times \omega$,  $D_{t,n}=\{\pp\in\Q:x_{t,n}\in\xbar_{\pp}\}$; and
\item  {\bf Henkin witnesses:}  For all $t\in\omega_1$, all $\<x_{s_i,n_i}:i<m\>$ with each $s_i\le t$ and all
$\phi(y,v_i:i<m)$,
$\{\pp\in\Q:\  \hbox{{\bf either}}\ \theta_{\pp}(\xbar_{\pp})\vdash\forall y\neg\phi(y,x_{s_i,n_i}:i<m) \ \hbox{{\bf or}}$ for some $n^*$,
$\theta_{\pp}(\xbar_{\pp})\vdash\phi(x_{t,n^*},x_{s_i,n_i}:i<m)\}$.
\item  For all $e\in M^*$, $D_e=\{\pp\in\Q:e\in \aabar_{\pp}$ and $\theta(\xbar_{\pp})\vdash x_{0,n}=e$ for some $n\in\omega\}$.
\end{enumerate}

\end{Lemma}

\bp  That each of these sets is open is immediate.  As for density, in all four clauses we will show that given some $\pp\in\Q$, we will find an extension $\qq\ge_{\Q}\pp$
with $\xbar_{\qq}$ a one-point extension of $\xbar_{\pp}$.  In all cases, we will put  $k_{\qq}:= k_{\pp}$, $\U_{\qq}=\U_{\pp}$ and since $\xbar_{\pp}$ is finite,
we can choose the color $c_{\qq}$ of the `new element' to be distinct from the other colors.  Because of that,  Clause~(8) for $\qq$ follows immediately from the fact $\pp\in\Q$.
Thus, for all four clauses, all of the work is in finding a striated type $\theta_{\qq}$ extending $\theta_{\pp}$.

(1)
 Fix $t\in\omega_1$ and choose
an arbitrary $\pp\in\Q$.  If $t\in u_{\pp}$ then there is nothing to prove, so assume otherwise.  Let $\ell=|u_{\pp}|$ and let $k=|\{s\in u_{\pp}:s<t\}|$.  Assume that $k<\ell$,
as the case of $k=\ell$ is similar, but easier.  Choose a perfect chain realization $(\Mbar,\bbar)$ of $\theta_{\pp}(\xbar_{\pp})$.  By Lemma~\ref{close}(2),
$\tp(\bbar_{<k}/\aabar_{\pp})$ is a striated type of length $k$.  By Lemma~\ref{large}(1), choose an $A^*$-large type $r\in S_{at}(\aabar_{\pp}\bbar_{<k})$ and choose a realization $e$ of $r$ in $M_k$.  One checks immediately that $\tp(\bbar_{<k} e/\aabar_{\pp})$ is a striated type of length $(k+1)$.  Now, also by Lemma~\ref{close}(2),
$(\Mbar_{\ge k},\bbar_{\ge k})$ is a perfect chain realization of $\tp(\bbar_{\ge k}/\aabar_{\pp}\bbar_{<k})$.  So, by Lemma~\ref{extend},
$(\Mbar_{\ge k},\bbar_{\ge k})$ is also a perfect chain realization of $\tp(\bbar_{\ge k}/\aabar_{\pp}\bbar_{<k}e)$.  In particular, $\tp(\bbar_{\ge k}/\aabar_{\pp}\bbar_{<k}e)$
is a striated type of length $(\ell-k)$.  Thus, by Lemma~\ref{close}(1), $\tp(\bbar_{<k}e\bbar_{\ge k}/\aabar_{\pp})$ is a striated type of length $(\ell+1)$.
Take $\aabar_{\qq}:=\aabar_{\pp}$, $\xbar_{\qq}:=\xbar_{\pp}\cup\{x_{t,0}\}$, 
and take $\theta_{\qq}(\xbar_{\qq})$ to be a complete formula in $\tp(\bbar_{<k}e\bbar_{\ge k}/\aabar_{\qq})$.

The proofs of (2) and (3) are extremely similar.  We prove (2) and indicate the adjustment necessary for (3).
Fix $(t,n)\in\omega_1\times\omega$.  By (1) and an inductive argument, we may assume we are given $\pp\in\Q$ with $t\in u_{\pp}$ and
$x_{t,n-1}\in\xbar_{\pp}$.  Say $|u_{\pp}|=\ell$ and assyne $t$ is the $(k-1)$st element of $u_p$ in ascending order.
Choose a perfect chain realization $(\Mbar,\bbar)$ of $\theta_{\pp}(\xbar_{\pp})$.  By Lemma~\ref{close}(2), $\tp(\bbar_{<k}/\aabar_{\pp})$ is striated of length $k$.
Choose an arbitrary $e\in M_k$\footnote{In the proof of (3), $e$ would be a realization of $\phi(y,b_{s_i,n_i}:i<m)$ in $M_k$, if one existed.} and adjoin it to $\bbar_{k-1}$.
More formally, let $\bbar_{<k}^*:=\<\bbar_j^*:j<k\>$, where $\bbar_j^*=\bbar_j$ for $j<k-2$, while $\bbar^*_{k-1}:=\bbar_{k-1}e$.
Note that $\tp(\bbar^*_{<k}/\aabar_{\pp})$ remains a striated type of length $k$.  
By Lemma~\ref{close}(2), $(\Mbar_{\ge k},\bbar_{\ge k})$ is a perfect chain realization of $\tp(\bbar_{\ge k}/\aabar_{\pp}\bbar_{<k})$.  So, by Lemma~\ref{extend}
it is also a perfect chain realization of $\tp(\bbar_{\ge k}/\aabar_{\pp}\bbar_{<k}^*)$.  In particular, $\tp(\bbar_{\ge k}/\aabar_{\pp}\bbar_{<k}^*)$ is a striated type of length 
$(\ell-k)$, so $\tp(\bbar^*_{<k}\bbar_{\ge k}/\aabar_{\pp})$ is a striated type of length $\ell$ extending $\theta_{\pp}(\xbar_{\pp})$.  Put $\xbar_{\qq}:=\xbar_{\pp}\cup\{x_{t,n}\}$
and let $\theta_{\qq}(\xbar_{\qq})$ be a complete formula isolating this type.

(4) is also similar and is left to the reader.
\endproof

The following Proposition follows immediately from the density conditions described above.

\begin{Proposition} \label{total} Let $G$ be a $\Q$-generic filter.  Then, in $V[G]$, a rich, $\U_G$-colored atomic model of $T$ exists,
where $\U_G=\{k\in\omega:k\in \U_{\pp}$ for some $\pp\in G\}$.
\end{Proposition}

\bp   There is a congruence $\sim_G$ defined on $X=\{x_{t,n}:t\in \omega_1,n\in\omega\}$ defined by $x_{t,n}\sim_G x_{s,m}$ if and only if
$\theta_{\pp}\vdash x_{t,n}=x_{s,m}$ for some $\pp\in G$.  Let $M_G$ be the model of $T$ with universe $X/\sim_G$ and relations 
$M_G\models\phi(a_1,\dots,a_k)$ if and only if there are $(x_{t_1,n_1},\dots,x_{t_k,n_k})\in X^k$ such that $[x_{t_i,n_i}]=a_i$ for each $i$
and  $\theta_{\pp}\vdash \phi(x_{t_1,n_1},\dots,x_{t_k,n_k})$ for some $\pp\in G$.  
Since $(\Q,\le_{\Q})$ has c.c.c., $M_G$ has size $\aleph_1$.  As notation, for each $t\in\omega_1$, let $M_{\le t}$ be the substructure of $M_G$
with universe $\{[x_{s,m}]:s\le t, m\in\omega\}$.  Then $M^*\preceq M_0$ and $M_{\le s}\preceq M_{\le t}\preceq M_G$ whenever $s\le t<\omega_1$.  
The definition of a striated type implies that $\tp([x_{t,0}]/A^*)$ is omitted in $M_{<t}$, hence 
the set $\{[x_{t,0}]:t\in\omega_1\}$ witnesses that $(M_G,\bbar^*)$ is
rich.  
Also, define $c_G:=\bigcup\{c_{\pp}:\pp\in G\}$.  Using the fact that each $\pp\in\Q$ is fully decided, check that $c_G$ is a $\U_G$-coloring of $(M_G,\bbar^*)$.
\endproof

Note that in the Conclusion below, such a $G\in V$ always exists, since $\B$ is countable.

\begin{Conclusion}  \label{quote}  Suppose $\B$ is a countable, transitive model of $ZFC^*$, with $\{M^*,T,L\}\subseteq \B$,
and let $G\in V$, $G\subseteq\Q$  be any filter meeting every dense $D\subseteq\Q$ with $D\in\B$.  Then:
  Let $\U_G=\{k\in\omega:k\in \U_{\pp}$ for some $\pp\in G\}$.  Then:
\begin{enumerate}
\item  $\U_G\in V$; and
\item  In $V$, there is a $\U_G$-colored, rich atomic model $(N,\bbar^*)$ of $T$.
\end{enumerate}
\end{Conclusion}

\bp  That $\U_G\in V$ is immediate, since both $\B$ and $G$ are.  As for (2), as $G$ meets every dense set in $\B$, $\B[G]$ is a countable, transitive model of $ZFC^*$, and
by applying Proposition~\ref{total},
$$\B[G]\models \hbox{`There is a rich, $\U_G$-colored $(M_G,\bbar^*)$ of size $\aleph_1$'}$$
Let $L'=L\cup\{c,R\}\cup\{c_m:m\in M^*\}$  Working in $\B[G]$, expand $M_G$ to an $L'$-structure $M'$, interpreting each $c_m$ by $m$,
interpreting  the unary function $c^{M'}$ as $c_G=\bigcup\{c_{\pp}:\pp\in G\}$, and the unary predicate $R^{M'}=\{[x_{t,0}]:t\in\omega_1\}$.

Now, for each $d,d'\in M'$ and $k\in\omega$, the relation $\tp_{M'}(d/A_k)=\tp_{M'}(d'/A_k)$ is definable by an $L'_{\omega_1,\omega}$-formula.
Thus, the binary function $\spl:(M')^2\rightarrow(\omega+1)$ is also $L'_{\omega_1,\omega}$-definable, hence, using the coloring $c$,
 there is an $L'_{\omega_1,\omega}$-sentence $\Psi$ 
stating that `$c$ induces a $\U_G$-coloring.'  Finally, using the $Q$-quantifier to state that $R$ is uncountable, there is an $L'_{\omega_1,\omega}$-sentence $\Phi\in\B[G]$
stating that the $L(\bbar^*)$-reduct of a given $L'$-structure is a rich, atomic model of $T$, that is $\U_G$-colored via $c$.
We finish by applying Proposition~\ref{trans} to $M'$ and $\Phi$.
\endproof

\subsection{Mass production}   \label{massp}
In this subsection we define a forcing $(\PP,\le_{\PP})$ such that a $\PP$-generic filter $G$ produces a perfect set $\{G_\eta:\eta\in 2^\omega\}$ of $\Q$-generic filters
such that the associated subsets $\{\U_{G_\eta}:\eta\in 2^\omega\}$ of $\omega$ are almost disjoint.
Although the application there is very different, the argument in this subsection is similar to one appearing in \cite{Sh480}.

We begin with one easy density argument concerning the partial $(\Q,\le_{\Q})$.
Fundamentally, it allows us to `stall' the construction for any fixed, finite length of time.

\begin{Lemma}  \label{mass}  For every  $\pp\in\Q$ and every $k^*>k_{\pp}$, there is $\qq\ge_{\Q}\pp$ such that $\xbar_{\qq}=\xbar_{\pp}$, 
(hence $c_{\qq}=c_{\pp}$); but  $k_{\qq}=k^*$ and $\U_{\qq}=\U_{\pp}$, i.e., $\U_{\qq}\cap [k_{\pp},k^*)=\emptyset$.
\end{Lemma}

\bp  Simply define $\qq$ as above and then verify that $\qq\in\Q$.
\endproof

\begin{Definition}  {\em  For $n\in\omega$, let 
$$\PP_n=\{(k,\pbar):k\in\omega, \pbar=\<p_\eta:\eta\in 2^n\>, \ \hbox{where each $p_\nu\in\Q$ and every $k_{p_\nu}=k$}\}$$
As notation, for $\pp\in\PP_n$, we let $k(\pp)$ denote the (integer) first coordinate of $\pp$.  For each $\ell<k(\pp)$, define the {\em trace of $\ell$,}
$\tr_\ell(\pp)=\{\nu\in 2^n:\ell\in \U_{p_{\nu}}\}$.

Let $\PP=\bigcup_{n\in\omega} \PP_n$.  As notation, for $\pp\in\PP$, $n(\pp)$ is the unique $n$ for which $\pp\in\PP_n$.
}
\end{Definition}

\begin{Definition}  {\em  Define an order $\le_{\PP}$ on $\PP$ by $\pp\le_{\PP} \qq$ if and only if
\begin{enumerate}
\item  $n(\pp)\le n(\qq)$, $k(\pp)\le k(\qq)$;
\item   $p_\nu\le_{\Q} q_\mu$ for all pairs $\nu\in 2^{n(\pp)}, \mu\in 2^{n(\qq)}$ satisfying $\nu\trianglelefteq\mu$; and
\item  For all $\ell\in [k(\pp),k(\qq))$, the set $\{\mu\mr{n(\pp)}:\mu\in \tr_\ell(\qq)\}$ is either empty or is a singleton.
\end{enumerate}
}
\end{Definition}

It is easily checked that $(\PP,\le_{\PP})$ is a partial order, hence a notion of forcing.  
The following Lemma describes the dense subsets of $\PP$.

\begin{Lemma} \label{PPdense}

\begin{enumerate}
\item  For each $n$ and $k$, $\{\pp\in\PP:n(\pp)\ge n\}$ and 
$\{\pp\in\PP: k(\pp)\ge k\}$ are dense;
\item   Suppose $D$ is a dense, open subset of $\Q$.  Then for every $n$ and every $\pp\in\PP_n$, there is $\qq\in\PP_n$
such that $\qq\ge_{\PP} \pp$ and, for every $\nu\in 2^n$, $\qq_\nu\in D$.
\end{enumerate}
\end{Lemma}

\bp Arguing by induction, it suffices to prove that for any given $\pp\in\PP$, there is $\qq\ge_\PP \pp$ with $n(\qq)=n(\pp)+1$ and
an $\rr\ge_{\PP}\pp$ with $k(\rr)>k(\pp)$.  Fix $\pp\in\PP$.  Say $\pp\in\PP_n$ and $\pp=(k,\pbar)$.  To construct $\qq$, for each $\nu\in 2^n$,
define $q_{\nu 0}=q_{\nu 1}=p_\nu$.  Let $\qbar:=\<q_\mu:\mu\in 2^{n+1}\>$ and $\qq=(k,\qbar)$.  Then $\qq\in\PP_{n+1}$ and $\qq\ge_{\PP} \pp$
(note that Clause~(3) in the definition of $\le_{\PP}$ is vacuously satisfied since $k(\pp)=k(\qq)$).  

To construct $\rr$, simply apply Lemma~\ref{mass} to each $p_\nu$ to produce an extension $r_\nu\ge_{\Q}p_\nu$ with $k_{r_\nu}=k+1$, but
$\U_{r_\nu}=\U_{p_\nu}$.  Then let $\rbar:=\<r_\nu:\nu\in 2^n\>$ and $\rr=(k+1,\rbar)$.  Then $\rr\ge_{\PP}\pp$ as required.

(2) Fix such a $D$ and $n$.  As we are working exclusively in $\PP_n$ and because $2^n$ is a fixed finite set, it suffices to prove that
for any chosen $\nu\in 2^n$, 
\begin{quotation}
\noindent For every $\pp\in\PP_n$ there is $\qq\in\PP_n$ with $\qq\ge_{\PP}\pp$ and $q_{\nu}\in D$.
\end{quotation}
To verify this, fix $\nu\in 2^n$ and $\pp\in \PP_n$.  Concentrating on $p_\nu$,  as $D$ is dense, choose $q_\nu\in D\cap\Q$ with $q_\nu\ge_{\Q}p_\nu$.
Let $k^*:=k_{q_\nu}$.  Next, for each $\delta\in 2^n$ with $\delta\neq \nu$, apply Lemma~\ref{mass} to $p_\delta$, obtaining some $q_\delta\in\Q$
satisfying $q_\delta\ge_{\Q} p_\delta$, $k_{q_\delta}=k^*$, but $\U_{q_\delta}=\U_{p_\delta}$.  Now, collect all of this data into a condition
$\qq\in \PP_n$ defined by $k(\qq)=k^*$ and $\qbar=\<q_\gamma:\gamma\in 2^n\>$, where each $q_\gamma$ is as above.  To see that $\qq\ge_{\PP}\pp$,
Clause~(3) is verified by noting that for every $\ell\in[k(\pp),k^*)$, $\tr_\ell(\qq)$ is either empty, or equals $\{\nu\}$, depending on whether or not $\ell\in \U_{q_\nu}$.
\endproof
 
 \begin{Notation}  \label{note}
 {\em  Suppose $\B\models ZFC^*$ and let $G^*\subseteq\PP$, $G^*\in V$
 be a filter meeting every dense subset $D^*\subseteq\PP$ with $D^*\in\B$.  
 For each $n$ and $\nu\in 2^n$, let
 $$G_\nu:=\{\pp\in\Q:\ \hbox{for some $\pp^*=(k,\pbar)\in G^*$, $\pp=\pp^*_\nu\}$}$$
 Then, for each $\eta\in 2^\omega$, let
 $$G_\eta:=\bigcup\{G_{\eta|n}:n\in\omega\}\ \ \hbox{and}\ \   \U_\eta:=\{\ell\in\omega:\ell\in \U_{\qq}\ \hbox{for some $\qq\in G_\eta$}\}$$
 }
 \end{Notation}

 \begin{Proposition}  \label{decompose}  In the notation of \ref{note}:
 \begin{enumerate}
 \item   For every $\eta\in 2^\omega$, $G_\eta\subseteq\Q$ is a filter meeting every dense $D\subseteq\Q$ with $D\in\B$; 
 \item  The sets $\{\U_\eta:\eta\in 2^\omega\}$ are an almost disjoint family of infinite subsets of $\omega$.
 \end{enumerate}
 \end{Proposition}
 
 \bp  (1) follows immediately from Lemma~\ref{PPdense}(2).
 
(2)  Choose distinct $\eta,\eta'\in 2^\omega$.  Choose $n_0$ such that $\eta|n\neq \eta'|n$ whenever $n\ge n_0$.  
By Lemma~\ref{PPdense}(1), choose $\pp^*\in G^*$ with $n(\pp^*)\ge n_0$. We show that $\U_\eta\cap \U_{\eta'}$ is finite by establishing that
if $\ell\in \U_{\eta}\cap\U_{\eta'}$, then $\ell\le k(\pp^*)$.

To establish this, choose $\ell\in \U_{\eta}\cap\U_{\eta'}$.  By unpacking the definitions, choose $\qq^*,\rr^*\in G^*$ such that, letting $\mu:=\eta |n(\qq^*)$
and $\mu':=\eta'| n(\rr^*)$, we have $\ell\in\U_{\qq^*_\mu}\cap\U_{\rr^*_{\mu'}}$.  As $G^*$ is a filter, choose $\ss^*\in G^*$ with $\ss^*\ge_{\PP} \pp^*,\qq^*,\rr^*$.
As notation, let $\delta:=\eta|n(\ss^*)$ and $\delta':=\eta'|n(\ss^*)$.

\medskip
\noindent
{{\bf Claim:}}  $\ell\in \U_{\ss^*_{\delta}}\cap \U_{\ss^*_{\delta'}}$.

\bp  As $\ell\in \U_{\qq^*_\mu}$, $\ell<k(\qq^*)$.  From $\qq^*\le_{\PP} \ss^*$ we conclude $k(\qq^*)\le k(\ss^*)$, so $\ell< k(\ss^*)$ as well.  
From $\qq^*\le_{\PP}\ss^*$ and $\mu\trianglelefteq \delta$ we obtain $\qq^*_\mu\le_{\Q} \ss^*_\delta$.  But then, as $\ell\in \U_{\qq^*_\mu}$, it follows
that $\ell\in \U_{\ss^*_{\delta}}$.  That $\ell\in \U_{\ss^*_{\delta'}}$ is analogous, using $\rr^*$ in place of $\qq^*$.

\medskip

Finally, assume by way of contradiction that $\ell\ge k(\pp^*)$.  The Claim implies that $\{\delta,\delta'\}\subseteq\tr_\ell(\ss^*)$.  As $\ell\in [k(\pp^*),k(\ss^*))$,
Clause~(3) of $\pp^*\le_{\PP} \ss^*$ implies that $\delta|n(\pp^*)=\delta'|n(\pp^*)$.  But, as $\eta|n(\pp^*)=\delta|n(\pp^*)$ and $\eta'|n(\pp^*)=\delta'|n(\pp^*)$, we contradict our choice of $\pp^*$.
\endproof

We close this section with the proof of Proposition~\ref{many}, which we restate for convenience.

\begin{Conclusion}  \label{family} There is a family $\{(N_\eta,\bbar^*):\eta\in 2^\omega\}$ of $2^{\aleph_0}$ rich, atomic models of $T$, each of size $\aleph_1$,
that are pairwise non-isomorphic over $\bbar^*$.
\end{Conclusion}

\bp
Choose any  countable, transitive model $\B$ of $ZFC^*$ and choose any $G^*\in V$, $G^*\subseteq\PP$, $G^*$ meets every dense subset $D^*\in \B$
(as $\B$ is countable, such a $G^*$ exists).    For each $\eta\in 2^\omega$, choose $G_\eta$ and $\U_\eta$ as in Proposition~\ref{decompose},
and apply Conclusion~\ref{quote} to get a rich $\U_\eta$-colored $(N_\eta,\bbar^*)$ in $V$.  That this family is pairwise non-isomorphic over $\bbar^*$
follows immediately from Corollary~\ref{criterion}, since the sets $\{\U_\eta:\eta\in 2^{\omega}\}$ are almost disjoint.
\endproof

\section{The proof of Theorem~\ref{big}}

Assume that the class $\At_T$ is not $\pcl$-small, as witnessed by an (uncountable) model $N^*$ containing a finite
tuple $\abar^*$.  Fix a countable, elementary substructure $M^*\preceq N^*$ that contains $\abar^*$.  To aid notation, let $D^*:=\pcl_{N^*}(\abar^*)$.  
We now split into cases, depending on the relationship between the cardinals $2^{\aleph_0}$ and $2^{\aleph_1}$.

\medskip
\noindent{\bf Case 1.}  $2^{\aleph_0}< 2^{\aleph_1}$.

\medskip
In this case, expand the language of $T$ to $L(D^*)$, adding a new constant symbol for each $d\in D^*$.   Then, the natural expansion $N^*_{D^*}$ $N^*$ to an 
$L(D^*)$-structure is a model of the infinitary $L(D^*)$-sentence $\Phi$ that entails $Th(N^*_{D^*})$ and ensures that every finite tuple is $L$-atomic with respect to $T$.
As $N^*_{D^*}$ is a model of $\Phi$ that realizes uncountably many types over the empty set (after fixing $D^*$!), it follows from \cite{Keislerbook}, Theorem~45 of  Keisler 
that there are $2^{\aleph_1}$ pairwise non-$L(D^*)$-isomorphic models $\Phi$, each of size $\aleph_1$.  As $2^{\aleph_0}< 2^{\aleph_1}$, it follows that there is a subfamily
of $2^{\aleph_1}$ pairwise non-$L$-isomorphic reducts to the original language $L$.  As each of these models are $L$-atomic,  we conclude that $\At_T$ has $2^{\aleph_1}$
non-isomorphic models of size $\aleph_1$.

\medskip\noindent{\bf Case 2.} $2^{\aleph_0}= 2^{\aleph_1}$.

\medskip  Choose $\bbar^*$ from $M^*$ as in Proposition~\ref{config} and apply Conclusion~\ref{family} to get a set $\F^*=\{(N_\eta,\bbar^*):\eta\in 2^\omega\}$ of atomic models,
each of size $\aleph_1$, that are pairwise non-isomorphic over $\bbar^*$.  Let $\F=\{N_\eta:\eta\in 2^\omega\}$ be the set of reducts of elements from $\F^*$.
By our cardinal hypothesis,  $\F$ has size $2^{\aleph_1}$.  The relation of $L$-isomorphism is an equivalence relation on $\F$, and
each $L$-isomorphism equivalence class has size at most $\aleph_1$ (since $\aleph_1^{<\omega}=\aleph_1$).  As $\aleph_1<2^{\aleph_1}$ we conclude that $\F$ has a 
subset of size $2^{\aleph_1}$ of pairwise non-isomorphic atomic models of $T$, each of size $\aleph_1$.
\endproof


\end{document}